\documentclass[11pt]{article}

\usepackage{amssymb}
\usepackage{amsmath}
\usepackage{theorem}
\usepackage{epsfig}
\usepackage{verbatim}
\usepackage{graphicx}

\newtheorem{thm}{Theorem}[section]

\newtheorem{prop}[thm]{Proposition}

\newcommand{\R}{\Bbb{R}}
\newcommand{\C}{\Bbb{C}}

\newcommand{\T}{\mathbb{T}}
\newcommand{\D}{\displaystyle}

\newcommand{\grad}{\nabla}

\newcommand{\al}{\alpha}

\newcommand{\g}{\gamma}

\newcommand{\da}{\partial_{\alpha}}
\newcommand{\pa}{\partial}
\begin{document}

\author{Angel Castro, Diego C\'ordoba, Charles Fefferman and Francisco Gancedo}
\title{Splash singularities for the one-phase\\
Muskat problem in stable regimes}

\date{}

\maketitle

\begin{abstract}
This paper shows finite time singularity formation for the Muskat problem in a stable regime. The framework we exhibit is with a dry region, where the density and the viscosity are set equal to $0$ (the gradient of the pressure is equal to $(0,0)$) in the complement of the fluid domain. The singularity is a splash-type: a smooth fluid boundary collapses due to two different particles evolve to collide at a single point. This is the first example of a splash singularity for a parabolic problem.
\end{abstract}

\maketitle


\section{Introduction}

This paper establishes some scenarios where the 2D Muskat problem produces splash singularities; that is to say, we prove that a free boundary evolving by the Muskat problem collapses at a single point while the interface prevails smooth. The situation is stable; we show geometries for initial data where the Rayleigh-Taylor condition holds.

The singularities we construct are ``splash'' singularities in which the interface self-intersects at a single point at the time of breakdown $T_{\ast}$ as in Fig. 1. Our previous papers \cite{ADCPJ} and \cite{ADCPJ2} showed the existence of a splash singularity for the water wave problem. The strategy there is to start with a ``splash" singularity at the time $T_{\ast}$ then solve water wave equation backwards in time. This yields a solution to the water wave equation in a time interval $\left[ T_{\ast}- \epsilon, T_{\ast}\right]$ that is well behaved at any time $\left[ T_{\ast}- \epsilon, T_{\ast}\right)$ but exhibits a splash at time $T_{\ast}$. In our present setting, we cannot use that strategy because the Muskat problem in the stable regime is parabolic and therefore cannot be solved backwards in time. The importance of this issue is made clear by the fact that water waves can form a ``splat" singularity \cite{ADCPJ2} whereas Muskat solution cannot \cite{CP}. (A ``splat" occurs when, at the time of breakdown, the interface self-intersects along an arc). On the other hand, an analysis of the Muskat problem has in common with our previous work on the water waves a conformal map to the ``tilde domain", see \cite{ADCPJ2}.

\begin{figure}[h]
\includegraphics[width=12cm,height=8cm]{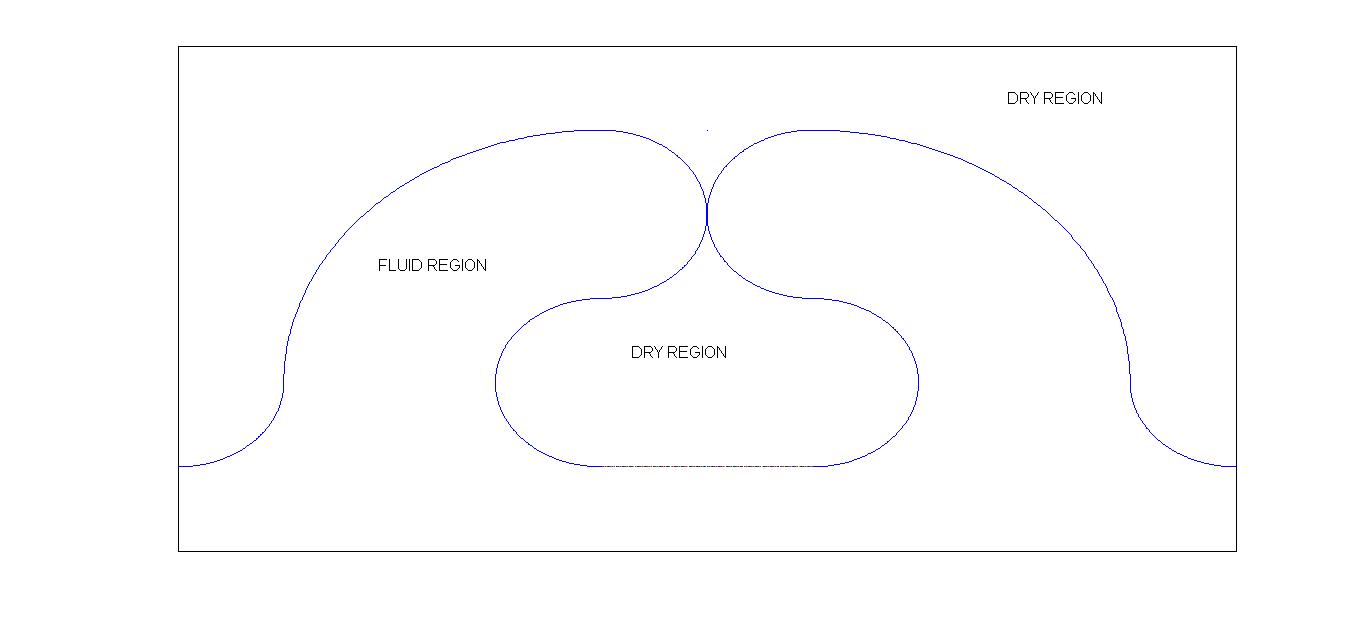}
\caption{Splash singularity}
\end{figure}

Recall the Muskat problem, which describes the evolution of two fluids of different nature in porous media. Both fluids are assumed to be immiscible and incompressible, been the most common example for applications the dynamics of water and oil \cite{Bear}. In two dimensions, the two fluids occupy the connected open set $D(t)$ and $\R^2\smallsetminus D(t)$ respectively. The characteristics of the fluids are their constant densities and viscosities. Then the step functions $\rho(x,t)$ and $\mu(x,t)$ represent the density and viscosity respectively on the media given by:
$$\rho(x,t)=\left\{\begin{array}{rl}
                  \rho_0,& x\in D(t),\\
                   \overline{\rho}_0,& x\in \R^2\smallsetminus D(t),
                \end{array}\right.
                $$
                $$
\mu(x,t)=\left\{\begin{array}{rl}
                  \mu_0,& x\in D(t),\\
                   \overline{\mu}_0,& x\in \R^2\smallsetminus D(t),
                \end{array}\right.
$$
for $x\in\R^2$, $t\geq 0$ and $\rho_0$, $\overline{\rho}_0$, $\mu_0$, $\overline{\mu}_0$ constant values. The main concern is about the dynamics of the common free boundary $\partial D(t)$, which is given by using the experimental Darcy's law:
\begin{equation}\label{DarcyIn}
\mu(x,t) v(x,t)=-\nabla p(x,t)-(0,\rho(x,t)).
\end{equation}
Here $v(x,t)=(v_1(x,t),v_2(x,t))$ is the incompressible velocity
\begin{equation}\label{IncomIn}
\nabla\cdot v(x,t)=0,
\end{equation}
and $p(x,t)$ is the scalar pressure. Above, the permeability of the media and the gravity constant are set equal to one without loss of generality.

The Muskat problem is a long standing matter \cite{Muskat} of recognized importance, specially because of its analogies with the evolutions of fluids in Hele-Shaw cells. In that setting the fluids are confined inside two closely parallel flat surfaces in such a way that the dynamics is essentially bidimensional. The Hele-Shaw evolution law is given by
\begin{equation*}
\frac{12}{b^2}\mu(x,t) v(x,t)=-\nabla p(x,t)-(0,\rho(x,t)),
\end{equation*}
where $b$ is distance among the surfaces. Therefore, it is possible to observe that for both different scenarios comparable phenomena and properties hold \cite{SCH}.

A main feature of the problem is the appearance of instabilities, which have been shown in different situations \cite{Otto},\cite{Lazlo}. From a contour dynamics point of view, the system of equations for the free boundary is essentially ill-possed from a Hadamard point of view \cite{SCH},\cite{DY}. Although taking into account surface tension effects the system becomes well-possed \cite{Esch}\cite{Ambrose3}, still it shows fingering \cite{Esch2} and exponential growing modes \cite{Guo}.

On the other hand, the Muskat problem is well-possed in stable regimes without surface tension \cite{DY},\cite{ADY},\cite{DRR}. This situation is reached for the problem when the difference of the gradient pressure jump at the free interface is positive \cite{Ambrose}. Then it is said that the Rayleigh-Taylor condition holds. In such case, linearizing the contour equation it is possible to get the following \cite{DY}:
\begin{equation}\label{lineal}
f^L_t(\al,t)=-\sigma\Lambda f^L(\al,t),
\end{equation}
where $(\al,f^L(\al,t))$ represents the free boundary ($\al\in\R$), $\sigma$ is the Rayleigh-Taylor constant and the operator $\Lambda$ is the square root of the negative Laplacian. Then, the fact that $\sigma>0$ turns the Muskat problem into a parabolic system at the linear level. This fact has been used to prove global in time regularity and instant analyticity for small initial data in different situations \cite{SCH},\cite{DY},\cite{Esch2},\cite{ccgs},\cite{Beck},\cite{Graner}.

For the case of equal viscosities ($\mu_0=\overline{\mu}_0$), the Rayleigh-Taylor condition holds when the more dense fluid lies below the interface and the less dense fluid lies above it \cite{DY}. In this situation, the regime is stable if the free boundary $\partial D$ is represented by the graph of a function $(\al,f(\al,t))$. In particular, it is possible to get a decay of the $L^{\infty}$ norm \cite{DP2} as follows:
$$
\|f\!-\!\frac{1}{2\pi}\int_{-\pi}^\pi f_0d\al\|_{L^\infty}(t)\leq \|f_0\!-\!\frac{1}{2\pi}\int_\T f_0d\al\|_{L^\infty}e^{-Ct},$$
for $f(\al+2\pi,t)=f(\al,t)$ and with $f(\al,t)\in L^2(\R)$
$$
\|f\|_{L^\infty}(t)\leq \|f_0\|_{L^\infty}(1+Ct)^{-1},\quad C=C(f_0)>0.
$$
It is easy to check that above formulas provide the same rate of decay than equation \eqref{lineal} for $f^L$ at the linear level. On the other hand, the $L^2$ norm evolution allows to control half derivative for $f^L$ due to the identity
$$
\|f^L\|^2_{L^2}(t)+2\sigma\int_0^t\|\Lambda^{1/2}f^L\|^2_{L^2}(s)ds=\|f_0^L\|^2_{L^2},
$$
meanwhile at the nonlinear level the following equality
$$
\|f\|^2_{L^2}(t)+\frac{\sigma}{\pi}\int_0^t\int_{\R}\int_{\R}
\ln \Big(1+\Big(\frac{f(\al,s)-f(\beta,s)}{\al-\beta}\Big)^2\Big)d\al d\beta ds
=\|f_0\|^2_{L^2},
$$
does not give a chance of gaining any regularity \cite{ccgs}.

The case of a drop on a solid substrate in porous media have been studied in \cite{Masmoudi}. This case considers the dynamics of one fluid, also known as the one-phase Muskat problem. The authors show local well-posedness on the problem with estimates independent of the contact angle.

In \cite{ADCPM} it is shown solutions of the Muskat equation for initial smooth stable graphs with precise geometries which enter in unstable regime becoming non-graph in finite time. The pattern is far from trivial and recently it has been shown to be richer for the inhomogeneous and confined problems (see \cite{JR} and references therein). In particular the significance of a turnover (non-graph scenario) is that the Rayleigh-Taylor condition breaks down. Furthermore, \cite{ADCP} there exist smooth initial data in the stable regime for the Muskat problem such that the solutions turn to the unstable regime and later the regularity breaks down. Therefore global-existence is false for some large initial data in the stable regime as the time evolution solutions develop singularities.

In this paper we show that the Muskat problem can develop singularities in stable regimes. The singularity is a splash, where for the free boundary given by
\begin{equation}\label{interdef}
\partial D(t)= \{z(\al,t)=(z_1(\al,t),z_2(\al,t)): \al\in\R\},
\end{equation}
there exist a blow-up time $T_s>0$ and a point $x_s\in\R^2$ such that $x_s=z(\alpha_1,T_s)=z(\alpha_2,T_s)$ for $\al_1\neq \al_2$. In particular the curve is regular, and satisfies the chord-arc condition up to the time $T_s$:
$$
|z(\al,t)-z(\beta,t)|\geq C_{ca}(t)|\al-\beta|,\quad \forall \al,\beta\in\R,\quad C_{ca}(t)>0,\quad t\in[0,T_s).
$$

 Free boundary incompressible fluid equations can develop splash singularities. This scenario have been shown for the incompressible Euler equations in the water waves form \cite{ADCPJ},\cite{ADCPJ2} which considers the evolution of a free
boundary given by air, with density 0, and water, with density 1, and irrotational velocity. This type of singularities can also be shown for the case with vorticity \cite{CS}. Although for the case of two incompressible fluids with positive densities, this setting has been recently ruled out \cite{FIL}. For Muskat this type of singularities does not also hold in the case in which $\mu_0=\overline{\mu}_0$ and $\rho_0\neq\overline{\rho}_0$ \cite{PB}.

In this work we show finite time splash singularities with $\overline{\rho}_0=\overline{\mu}_0=0$:
\begin{equation}\label{rmIn}
(\rho(x,t),\mu(x,t))=\left\{\begin{array}{rl}
                  (\rho_0,\mu_0)& x\in D(t),\\
                   (0,0),& x\in \R^2\smallsetminus D(t),
                \end{array}\right.
\end{equation}
considering the one fluid dynamics with $\R^2\smallsetminus D(t)$ a dry region. We also yield some geometries for the interface where the Rayleigh-Taylor condition is satisfied, getting rid of unstable situations. The main theorem of the paper is the following:

\begin{thm}\label{splashvacuum}
There exist an open set of curves $\mathcal{O}\subset H^3$, satisfying the chord-arc and Rayleigh-Taylor condition, such that for any $z_0\in\mathcal{O}$ the solution of Muskat (\ref{DarcyIn},\ref{IncomIn},\ref{interdef},\ref{rmIn}) with $z(\al,0)=z_0(\al)$ violates the chord-arc condition at a finite time $T_s=T_s(z_0)>0$. In addition, this holds in such a way that $z(\al_1,T_s)=z(\al_2,T_s)$ with $\al_1\neq\al_2$. At the time $T_s$ the Muskat system (\ref{DarcyIn},\ref{IncomIn},\ref{interdef},\ref{rmIn}) breaks down.
\end{thm}

In the rest of the paper we show the proof of above result splitting it in several sections. In section 2 we construct a family of curves $z^l$ for which there is a unique self-intersection point $x_s$ where $x_s=z^l(\al_1)=z^l(\al_2)$ with $\al_1\neq\al_2$ and $\pa_\alpha z^l_1(\alpha_1)=\pa_\alpha z^l_1(\alpha_2)=0$. Plugging these curves in Darcy's law, we get that the Rayleigh-Taylor condition holds. Furthermore, the velocity indicates that the self-intersection point is going to disappear going backward in time. A more general scenario can be found in section \ref{generales}. In section 3 we show how to make sense the problem with a self-intersecting point, transforming the Muskat problem into a new contour dynamics equation we call $P$(Muskat). Up to the time of the splash we can recover Muskat from $P$(Muskat), but at the time of splash $P$(Muskat) makes sense and it is possible to go further in time. In section 4 we prove local existence of the $P$(Muskat) system. In section
 5 we show a stability result for $P$(Muskat). Finally, in section 6 we show how the family of curves $z^l(\al)$ together with the existence and stability for $P$(Muskat) allow us to conclude the proof of theorem \ref{splashvacuum}.

\section{Self-intersecting stable curves with suitable sign of velocity}\label{initialdata}

In this section we show that there exits a family of splash curves such that the Rayleigh-Taylor condition hold and with velocities which separate the splash point running backward-in-time.

First we use Hopf's lemma to achieve the Rayleigh-Taylor condition. Taking divergence in Darcy's law \eqref{DarcyIn} we have
$$\Delta p(x,t)=0,$$
for any $x\in D(t)$.
In addition, the continuity of the pressure on the free boundary \cite{ADY} and the fact that
$$
-\grad p(x,t)=(0,0)
$$
for any $x$ in the interior of $\R^2\smallsetminus D(t)$ allow us to get
 $$p(z(\al,t),t)=0.$$
On the other hand we consider velocities with mean zero vorticity $\partial_{x_1}v_2-\partial_{x_2}v_1$ which provides $v\in L^2(D(t))$ and finite energy settings. Approaching to infinity in $D(t)$ yields
\begin{align*}
\lim_{x_2\to -\infty} v(x,t)=0,
\end{align*}
and therefore Darcy's law gives
\begin{align*}
\lim_{x_2\to-\infty}\partial_{x_1} p(x,t)=0,\\
\lim_{x_2\to-\infty}\partial_{x_2} p(x,t)=-\rho_0.
\end{align*}
It is possible to find that $p(x,y)\sim -\rho_0 x_2+c(t)$ when $x_2\to -\infty$ and to conclude that the pressure is positive in $D(t)$ by the maximum principle for harmonic functions. In this situation we can apply Hopf's lemma to obtain that
\begin{align}\label{Hopf}
-\grad p(z(\al,t),t)\cdot\partial_\al^\bot z(\al,t) \geq c(t)>0,
\end{align}
where $\partial_\al^\bot z(\al,t)=(-\partial_\al z_2(\al,t),\partial_\al z_1(\al,t))$ is the normal vector pointing out the domain $D(t)$.

Next we deal with curves $z^l(\al)$ with a splash point $x_s=z^l(\al_1)=z^l(\al_2)$
for $\al_1\neq\al_2$ where
$$\partial_\alpha z^l_1(\alpha_1)=\partial_\alpha z^l_1(\alpha_2)=0.$$
We show that this configuration provides a sign for the velocity at $x_s$. Taking the trace of Darcy's law to the surface and multiplying by $\partial_\al^\bot z^l(\al)$ we have that
\begin{align*}
\mu_0 v(z^l(\al))\cdot \partial_\al^\bot z^l(\al)= -\grad p(z^l(\alpha))\cdot \partial_\al^\bot z^l(\al) -\rho_0\partial_\alpha z_1^l(\alpha).
\end{align*}
Thanks to our choice of the splash curve it must be satisfied
\begin{align}\label{last}
v(z^l(\al_i))\cdot \partial_\al^\bot z^l(\al_i)= -\mu_0^{-1}\grad p(z^l(\alpha_i))\cdot \partial_\al^\bot z^l(\al_i)\geq c>0,\quad i=1,2,
\end{align}
where again we have used Hopf's lemma \eqref{Hopf}. It is clear that \eqref{last} implies that the velocity separates the splash point backwards in time.

In figure 1 we give a graphic sketch of the kind of splash singularities we are considering.


Theses splash curves yield the simplest scenario we can consider. In section \ref{generales} we show the existence of different geometries that give rise to a splash singularity for the one-phase Muskat problem.

\section{Transformation to a non-splash scenario}

This section is devoted to transform the system into a new contour evolution equation where we handle the splash singularity. We consider solutions of Muskat satisfying (\ref{DarcyIn},\ref{IncomIn},\ref{interdef},\ref{rmIn}) for regular $z(\al,t)$ satisfying the chord-arc condition. Taking limit as $x\to z(\al,t)$ from $D(t)$ we find
$$
v(z(\al,t),t)=u(\al,t),
$$
where
$$
u(\al,t)=BR(z,\omega)(\al,t)+\frac{\omega(\al,t)}{2}\frac{z_\al(\al,t)}{|z_\al(\al,t)|^2}.
$$
$BR$ stands for the Birkhoff-Rott integral, which is given by
\begin{align}
\begin{split}\label{BR}
BR(\al,t)=BR(z,\omega)(\al,t)&=\frac{1}{2\pi}PV\int_{\R}\frac{(z(\al,t)-z(\al-\beta,t))^{\bot}}{|z(\al,t)-z(\al-\beta,t)|^2}\omega(\al-\beta,t)d\beta,
\end{split}
\end{align}
and $\omega$ is the amplitude of the vorticity concentrated on the free boundary:
$$
(\partial_{x_1}v_2-\partial_{x_2}v_1)(x,t)=\omega(\beta,t)\delta(x=z(\beta,t)).
$$
By approaching to the contour in Darcy's law and taking the dot product with $z_\al(\al,t)$ it is easy to relate the amplitude of the vorticity and the free boundary by an elliptic implicit equation:
\begin{align}
\begin{split}\label{fw}
\omega(\al,t)&=-2BR(z,\omega)(\al,t)\cdot
\da z(\al,t)-2\frac{\rho_0}{\mu_0}\da
z_2(\al,t).
\end{split}
\end{align}
We have the dynamics given by the following contour equation
\begin{align}
\begin{split}\label{Muskat}
\D z_t(\al,t)&=u(\al,t)+c(\al,t)\da z(\al,t)
\end{split}
\end{align}
where $c$ represents reparameterization freedom. See \cite{ADY} for a detail derivation of the system.

In a periodic setting in the $x_1$ direction, we will transform the system with the conformal map:
$$
P(w)=\big(\tan (w/2)\big)^{1/2},\quad w\in\C.
$$
Above, the branch of the square root is chosen in such a way that crosses the self-intersecting point of the $z^l(\al)$ curve from before section. Therefore $P(z^l(\al))$ becomes a one-to-one curve.

 We then consider by this new transformation the curve $\tilde{z}(\al,t)=P(z(\al,t))$. This provides easily
$$
\tilde{z}_\al(\al,t)=\nabla P(z(\al,t))z_\al(\al,t),
$$
and
$$
\tilde{z}_t(\al,t)=\nabla P(z(\al,t))z_t(\al,t)=\nabla P(z(\al,t))(u(\al,t)+c(\al,t)z_\al(\al,t))=
$$
$$
=\nabla P(z(\al,t))u(\al,t)+c(\al,t)\tilde{z}_\al(\al,t).
$$
For the potential $\phi(x,t)$ ($\grad \phi(x,t) = v(x,t)$) we define in the tilde domain $\tilde{\phi}(\tilde{x},t)=\phi(x,t)$.
Then
$$
v(x,t)=\grad\phi(x,t)=(\grad\tilde{\phi})(P(x),t)\grad P(x)=\grad P(x)^T(\grad\tilde{\phi})(P(x),t).
$$
Taking limit we find
$$
u(\al,t)=\grad P(z(\al,t))^T(\grad\tilde{\phi})(P(z(\al,t)),t)=\grad P(z(\al,t))^T\tilde{u}(\al,t),
$$
where $\tilde{u}(\al,t)=\grad\tilde{\phi}(\tilde{z}(\al,t),t)$. It yields
\begin{equation}\label{ecu}
\tilde{z}_t(\al,t)=Q^2(\al,t)\tilde{u}(\al,t)+c(\al,t)\tilde{z}_\al(\al,t),
\end{equation}
where $Q^2$ is given by
$$
\grad P(z(\al,t))\grad P(z(\al,t))^T=Q^2(\al,t) I,
$$
and $I$ is the $2\times2$ identity matrix. In other words
\begin{equation}\label{Q}
Q^2(\al,t)=\Big|\frac{dP}{dw}(z(\al,t))\Big|^2=\Big|\frac{dP}{dw}(P^{-1}(\tilde{z}(\al,t)))\Big|^2.
\end{equation}

Next we consider the velocity $\tilde{v}$ defined on the whole space by
$$
\tilde{v}(\tilde{x},t)=\grad\tilde{\phi}(\tilde{x},t)=
\frac{1}{2\pi}PV\int_{\R}\frac{(\tilde{x}-\tilde{z}(\al-\beta,t))^{\bot}}{|\tilde{x}-\tilde{z}(\al-\beta,t)|^2}\tilde{\omega}(\al-\beta,t)d\beta,
$$
where
$$
(\partial_{\tilde{x}_1}\tilde{v}_2-\partial_{\tilde{x}_2}\tilde{v}_1)(\tilde{x},t)=\tilde{\omega}(\beta,t)\delta(\tilde{x}=\tilde{z}(\beta,t)),
$$
in a distributional sense. Approaching to the free boundary it is possible to obtain
\begin{equation}\label{utilde}
\tilde{u}=BR(\tilde{z},\tilde{\omega})+\frac{\tilde{\omega}}{2|\tilde{z}_{\al}|^2}\tilde{z}_{\al}.
\end{equation}
In order to close the system we integrate Darcy's law to find
$$
\mu_0 \phi(z(\al,t),t)=-p(z(\al,t),t)-\rho_0 z_2(\al,t)=-\rho_0 z_2(\al,t),
$$
due to the continuity of the pressure at the free boundary and the vacuum state. The conformal map $P$ provides
\begin{equation}\label{DarcyP}
\mu_0\tilde{\phi}(\tilde{z}(\al,t),t)=-\rho_0 P_2^{-1}(\tilde{z}(\al,t)).
\end{equation}
Taking one derivative and identity \eqref{utilde} allow us to find
$$
\mu_0(BR(\tilde{z},\tilde{\omega})\cdot\tilde{z}_{\al}+\frac{\tilde{\omega}}{2})=-\rho_0\da(P_2^{-1}(\tilde{z})).
$$
We rewrite above identity as
\begin{equation}\label{omegatilde}
\tilde{\omega}(\al,t)=-2 BR(\tilde{z},\tilde{\omega})(\al,t)\cdot\tilde{z}_{\al}(\al,t)-2\frac{\rho_0}{\mu_0}\da(P_2^{-1}(\tilde{z}(\al,t))).
\end{equation}
Next we will pick a tangential component to get $|\tilde{z}_\al|$ depending only on the variable $t$. Identities \eqref{ecu} and \eqref{utilde} give
\begin{equation}\label{ecBR}
\tilde{z}_t(\al,t)=Q^2(\al,t)BR(\tilde{z},\tilde{\omega})(\al,t)+\tilde{c}(\al,t)\tilde{z}_\al(\al,t),
\end{equation}
for $\tilde{c}=Q^2\tilde{\omega}/(2|\tilde{z}_{\al}|^2)+c$. This provides
\begin{equation}\label{ctilde}
\tilde{c}(\al,t)=\frac{\al+\pi}{2\pi}\int_{-\pi}^\pi \partial_\beta(Q^2BR)(\beta,t)\cdot\frac{\tilde{z}_\beta(\beta,t)}{|\tilde{z}_\beta(\beta,t)|^2}d\beta-
\int_{-\pi}^\alpha \partial_\beta(Q^2BR)(\beta,t)\cdot\frac{\tilde{z}_\beta(\beta,t)}{|\tilde{z}_\beta(\beta,t)|^2}d\beta.
\end{equation}
We end up with a contour equation given by (\ref{omegatilde},\ref{ecBR},\ref{ctilde}).

Finally we will find the Rayleigh-Taylor condition in terms of $\tilde{z}$. We define $\tilde{p}(\tilde{x},t)=p(x,t)$ to obtain with Darcy's law
$$
-\grad\tilde{p}(\tilde{x},t)=\mu_0\grad\tilde{\phi}(\tilde{x},t)+\rho_0 \grad P_2^{-1}(\tilde{x}).
$$
Approaching to the free boundary, we find easily
\begin{equation}\label{RT}
\tilde{\sigma}(\al,t)=-\grad\tilde{p}(\tilde{z}(\al,t),t)\cdot\tilde{z}_{\al}^\bot=\mu_0 BR(\tilde{z},\tilde{\omega})\cdot\tilde{z}_{\al}^\bot+\rho_0 \grad P_2^{-1}(\tilde{z}(\al,t))\cdot\tilde{z}_{\al}^\bot.
\end{equation}

\section{Local-existence in the tilde domain}

This section is devoted to prove local existence for $\tilde{z}$ solutions of (\ref{omegatilde},\ref{ecBR},\ref{ctilde}) with $\tilde{z}\in C([0,T];H^k)$ with $k\geq 3$.
We show the proof for $k=3$ with the rest of the cases being analogous. In order to simplified the exposition  we suppress the time variable and the tilde in the equation. We follow the same strategy as in \cite{ADY}. We define
$$
q^0=(0,0),\quad q^1=(\frac{1}{\sqrt{2}},\frac{1}{\sqrt{2}}),\quad q^2=(\frac{-1}{\sqrt{2}},\frac{1}{\sqrt{2}}),\quad q^3=(\frac{-1}{\sqrt{2}},\frac{-1}{\sqrt{2}}),\quad q^4=(\frac{1}{\sqrt{2}},\frac{-1}{\sqrt{2}}),
$$
which are the singular points of the $P^{-1}$ conformal map. We set $z(\al,t)$ to hold $\tilde{z}(\al,t)\neq q^l$ for $l=0,...,4$. In order to get this we fix $\overline{D(0)}$ so that $\frac{dP}{dw}(w)\neq 0$ for any $w\in \overline{D(0)}$ without loss of generality. We will check that this property remains true for short time. Next we define the quantity
\begin{equation}\label{E}
E_k(z,t)=E_k(t)=\|z\|^2_{H^k}(t)+\|F(z)\|_{L^\infty}^2(t)+\frac{1}{m(Q^2\sigma)(t)}+\sum_{l=0}^4\frac{1}{m(q^l)(t)},
\end{equation}
where
$$
F(z)=\frac{|\beta|}{|z(\al)-z(\al-\beta)|},\quad \al,\beta\in[-\pi,\pi],
$$
and
$$
m(Q^2\sigma)(t)=\min_{\al\in\T}Q^2(\al,t)\sigma(\al,t),\quad m(q^l)(t)=\min_{\al\in\T}|z(\al,t)-q^l|.
$$
We shall show a proof of the following result:
\begin{prop}
Let $z(\al,t)$ be a solution of (\ref{omegatilde},\ref{ecBR},\ref{ctilde}). Then, the following estimate holds:
$$
\frac{d}{dt}E_k(t)\leq C (E_k(t))^p
$$
for $k\geq 3$. The constants $C$ and $p$ depend only on $k$.
\end{prop}

Below we will show the proof for $k=3$, being the rest of the cases analogous. These a priori estimates will lead to a local existence result for the contour equation in the tilde domain.

We refer the reader to \cite{ADY} in order to obtain
$$
\frac{d}{dt}\Big(\|z\|^2_{L^2}(t)+\|F(z)\|_{L^\infty}^2(t)+\frac{1}{m(Q^2\sigma)(t)}+\sum_{l=0}^4\frac{1}{m(q^l)(t)}\Big)\leq  C (E_k(t))^p,
$$
as a similar approach can be made. Next we check
$$
\frac{d}{dt}\|\da^3z\|^2_{L^2}(t)=2\int \da^3 z(\al)\cdot\da^3 z_t(\al)d\al.
$$
We can estimate most of the terms as in \cite{ADY}. We also quote \cite{ADCPJ} for dealing with the $Q^2$ factor. This term do not introduce any unbounded character as
$$
\|Q^2\|_{H^k}\leq C(E_k(t))^p.
$$
We will show how to deal with the unbounded and therefore singular terms. We find
$$
\frac{d}{dt}\|\da^3z\|^2_{L^2}(t)\leq C (E_k(t))^p+I,
$$
for
$$
I=\int \da^3 z(\al)\cdot Q^2(\al)\frac{1}{\pi}\int \frac{(z(\al)-z(\al-\beta))^{\bot}}{|z(\al)-z(\al-\beta)|^2}\da^3\omega(\al-\beta)d\beta d\al.
$$
We get $I\leq C (E_k(t))^p+II$ where
$$
II=\int \da^3 z(\al)\cdot \frac{z_\al^\bot(\al)}{|z_\al(\al)|^2}Q^2(\al)H(\da^3\omega)(\al)d\al.
$$
Identity $H(\da)=\Lambda$ allows us to rewrite $II$ as follows
$$
II=\frac{1}{|z_\al(\al)|^2}\int  \Lambda(\da^3z\cdot z_\al^\bot Q^2)(\al)\da^2\omega(\al)d\al.
$$
Next we can use formula \eqref{omegatilde} to split further $II=III+IV$ where
$$
III=\frac{-2}{|z_\al(\al)|^2}\int  \Lambda(\da^3z\cdot z_\al^\bot Q^2)(\al)\da^2(BR(z,\omega)\cdot z_\al)(\al)d\al,
$$
and
$$
IV=\frac{-2\rho_0\mu_0^{-1}}{|z_\al(\al)|^2}\int  \Lambda(\da^3z\cdot z_\al^\bot Q^2)(\al)\da^3(P_2^{-1}(z))(\al)d\al.
$$
The term $III$ can be estimated as $K_3$ in pg. 514 of \cite{ADY}. An analogous approach provides
\begin{equation}\label{rt1}
III\leq  C (E_k(t))^p-\frac{2}{|z_\al(\al)|^2}\int Q^2(\al)BR(z,\omega)(\al)\cdot z_\al^\bot(\al) \da^3 z(\al)\cdot\Lambda(\da^3z)(\al)d\al.
\end{equation}
For $IV$ we consider the most singular terms as the rest are bounded: $IV\leq C (E_k(t))^p+V$ where
$$
V=-\frac{2\rho_0\mu_0^{-1}}{|z_\al(\al)|^2}\int  \Lambda(\da^3z\cdot z_\al^\bot Q^2)(\al)(\grad P_2^{-1})(z(\al))\cdot\da^3z(\al)d\al.
$$
Then we split further $V=VI+VII+VIII+IX$ by writing the components of the curve:
$$
VI=\frac{2\rho_0\mu_0^{-1}}{|z_\al(\al)|^2}\int  \Lambda(\da^3z_1\da z_2 Q^2)(\al)\partial_{\tilde{x}_1}P_2^{-1}(z(\al))\da^3z_1(\al)d\al,
$$
$$
VII=\frac{2\rho_0\mu_0^{-1}}{|z_\al(\al)|^2}\int  \Lambda(\da^3z_1\da z_2 Q^2)(\al)\partial_{\tilde{x}_2}P_2^{-1}(z(\al))\da^3z_2(\al)d\al,
$$
$$
VIII=-\frac{2\rho_0\mu_0^{-1}}{|z_\al(\al)|^2}\int  \Lambda(\da^3z_2 \da z_1 Q^2)(\al)\partial_{\tilde{x}_1}P_2^{-1}(z(\al))\da^3z_1(\al)d\al,
$$
$$
IX=-\frac{2\rho_0\mu_0^{-1}}{|z_\al(\al)|^2}\int  \Lambda(\da^3z_2 \da z_1 Q^2)(\al)\partial_{\tilde{x}_2}P_2^{-1}(z(\al))\da^3z_2(\al)d\al.
$$
The commutator estimate
$$
\|\Lambda(gf)-g\Lambda f\|_{L^2}\leq C\|g\|_{C^{1,\frac13}}\|f\|_{L^2},
$$
yields
\begin{equation}\label{rt21}
VI\leq C (E_k(t))^p+\frac{2\rho_0\mu_0^{-1}}{|z_\al(\al)|^2}\int Q^2(\al) \partial_{\tilde{x}_1}P_2^{-1}(z(\al)) \da z_2(\al) \da^3z_1(\al)\Lambda(\da^3z_1)(\al)d\al,
\end{equation}
and
\begin{equation}\label{rt22}
IX\leq C (E_k(t))^p-\frac{2\rho_0\mu_0^{-1}}{|z_\al(\al)|^2}\int Q^2(\al) \partial_{\tilde{x}_2}P_2^{-1}(z(\al)) \da z_1(\al)\da^3z_2(\al)\Lambda(\da^3z_2)(\al)d\al.
\end{equation}
Similarly for $VII$:
$$
VII\leq C (E_k(t))^p+\frac{2\rho_0\mu_0^{-1}}{|z_\al(\al)|^2}\int Q^2(\al)\partial_{\tilde{x}_2}P_2^{-1}(z(\al)) \da z_2(\al) \da^3z_2(\al)\Lambda(\da^3z_1)(\al)d\al.
$$
Identity
$$
\da z_2(\al) \da^3z_2(\al)=-\da z_1(\al) \da^3z_1(\al)+|\da^2z(\al)|^2,
$$
provides
\begin{equation}\label{rt23}
VII\leq C (E_k(t))^p-\frac{2\rho_0\mu_0^{-1}}{|z_\al(\al)|^2}\int Q^2(\al)\partial_{\tilde{x}_2}P_2^{-1}(z(\al)) \da z_1(\al) \da^3z_1(\al)\Lambda(\da^3z_1)(\al)d\al.
\end{equation}
Proceeding in a similar manner we can get
\begin{equation}\label{rt24}
VIII\leq C (E_k(t))^p+\frac{2\rho_0\mu_0^{-1}}{|z_\al(\al)|^2}\int Q^2(\al)\partial_{\tilde{x}_1}P_2^{-1}(z(\al)) \da z_2(\al) \da^3z_2(\al)\Lambda(\da^3z_2)(\al)d\al.
\end{equation}
Adding the inequalities \eqref{rt21}, \eqref{rt22}, \eqref{rt23} and \eqref{rt24} it is easy to get
\begin{equation*}\label{rt2}
V\leq C (E_k(t))^p-\frac{2\rho_0\mu_0^{-1}}{|z_\al(\al)|^2}\int Q^2(\al)\grad P_2^{-1}(z(\al))\cdot z_\al^\bot(\al) \da^3z(\al)\cdot\Lambda(\da^3z)(\al)d\al.
\end{equation*}
Above inequality together with \eqref{rt1} let us obtain
$$
III\leq C (E_k(t))^p-\frac{2\mu_0^{-1}}{|z_\al(\al)|^2}\int Q^2(\al)\sigma(\al)\da^3z(\al)\cdot\Lambda(\da^3z)(\al)d\al
$$
with $\sigma$ given in \eqref{RT}.

Finally we obtain
$$
\frac{d}{dt}\|\da^3z\|^2_{L^2}(t)\leq C (E_k(t))^p-\frac{2\mu_0^{-1}}{|z_\al(\al)|^2}\int Q^2(\al)\sigma(\al)\da^3z(\al)\cdot\Lambda(\da^3z)(\al)d\al.
$$
From the a priori energy estimates we have that $m(Q^2\sigma)(t)>0$ which together with the pointwise inequality
$
2f\Lambda (f)\geq \Lambda(f^2)
$
yields
$$
\frac{d}{dt}\|\da^3z\|^2_{L^2}(t)\leq C (E_k(t))^p-\frac{\mu_0^{-1}}{|z_\al(\al)|^2}\int Q^2(\al)\sigma(\al)\Lambda(|\da^3z|^2)(\al)d\al.
$$
Integration by parts for the $\Lambda$ operator gives the desired estimate.

\section{Stability for the Muskat problem}

This section is devoted to show the proof of the following result:

\begin{prop}
Let $x(\al,t)$ and $y(\al,t)$ be two curves which satisfy the contour equation (\ref{omegatilde},\ref{ecBR},\ref{ctilde}). Then, the following estimate holds:
$$
\frac{d}{dt}\|x-y\|_{H^1}(t)\leq C (\sup_{[0,T]}E_3(x,t)+\sup_{[0,T]}E_3(y,t))^p\|x-y\|_{H^1}(t).
$$
Above $E_3(x,t)$ and $E_3(y,t)$ are given by \eqref{E}. The constants $C$ and $p$ are universal.
\end{prop}

\emph{Proof:}
In order to simplified the exposition we suppress the time variable and we denote $f'=f(\al-\beta)$, $f=f(\al)$, $f_-=f-f'$ and $\int=\int_{\T}$.

We consider two solutions of the system $x(\al,t)$ and $y(\al,t)$ in $C([0,T];H^3(\T))$ with $\gamma$ and $\zeta$ its vorticity amplitudes given by \eqref{omegatilde}. We will also denote by $Q^2_x$, $Q^2_y$, $BR_x$, $BR_y$ and $c_x$, $c_y$ the factors $Q^2$, Birhoff-Rott integrals and parametrization constants associated to $x$ and $y$ respectively (see \eqref{Q}, \eqref{BR} and \eqref{ctilde}). During the time of existence $T>0$ one finds $\sup_{[0,T]}E_3(x,t)$ and $\sup_{[0,T]}E_3(y,t)$ bounded so that we will write
$$
C (\sup_{[0,T]}E_3(x,t)+\sup_{[0,T]}E_3(y,t))^p\leq C
$$
by abuse of notation.

 For the function $z(\al,t)=x(\al,t)-y(\al,t)$ one finds
$$
\frac{1}{2}\frac{d}{dt}\|z\|^2_{L^2}=\int z\cdot z_t d\al= I_1+I_2+I_3+I_4,
$$
where
$$
I_1=\int z\cdot (Q^2_x-Q^2_y)BR_xd\al,\quad
I_2=\int z\cdot Q^2_y(BR_x-BR_y)d\al,
$$
$$
I_3=\int z\cdot (c_x-c_y)x_\al d\al,\quad
I_4=\int z\cdot c_y z_\al d\al.
$$
Then for $I_1$ we find
$$
I_{1}\leq \|z\|_{L^\infty}\|Q^2_x-Q^2_y\|_{L^2}\|BR_x\|_{L^2}\leq C\|z\|^2_{H^1}.
$$
In $I_2$ we split further as follows:
$$
I_{2,1}=\frac1{2\pi}\int z\cdot Q^2_y\int\frac{z^\bot_-}{|x_-|^2}\gamma'd\beta d\al,\quad
I_{2,2}=\frac1{2\pi}\int z\cdot Q^2_y\int y^\bot_-(\frac{1}{|x_-|^2}-\frac{1}{|y_-|^2})\gamma'd\beta d\al,
$$
$$
I_{2,3}=\int z\cdot Q^2_y BR(y,\omega)d\al,
$$
where $\omega=\gamma-\zeta$. In $I_{2,1}$, for the integral in $\beta$, we find a kernel of degree $-2$ applied to $z$ thus
$$
I_{2,1}\leq C\|z\|_{H^1}^2.
$$
Since
$$
I_{2,2}=\frac{-1}{2\pi}\int z\cdot Q^2_y\int y^\bot_-\frac{(x_-+y_-)\cdot z_-}{|x_-|^2|y_-|^2}\gamma'd\beta d\al,
$$
we recognize again a kernel of degree $-2$ applied to $z$ above so that
$$
I_{2,2}\leq C\|z\|_{H^1}^2.
$$
For $I_{2,3}$ it is easy to check that $BR$ has a kernel of degree $-1$ and therefore
$$
I_{2,3}\leq C\|z\|_{L^2}\|\omega\|_{L^2}.
$$
In order to deal with $\|\omega\|_{L^2}$ we write
$$
\omega+2BR(x,\omega)\cdot x_\al=2BR(y,\zeta)\cdot y_\al-2BR(x,\zeta)\cdot x_\al+2\frac{\rho_0}{\mu_0}(\grad P_2^{-1}(y)\cdot y_{\al}-\grad P_2^{-1}(x)\cdot x_{\al}).
$$
Bounds for the operator $(I+2BR(x,\cdot)\cdot x_\al)^{-1}$ (see \cite{ADY}) allow us to get
$$
\|\omega\|_{L^2}\leq C\|2BR(y,\zeta)\cdot y_\al-2BR(x,\zeta)\cdot x_\al+\frac{\rho_0}{\mu_0}(\grad P_2^{-1}(y)\cdot y_{\al}-\grad P_2^{-1}(x)\cdot x_{\al})\|_{L^2}.
$$
We proceed as before to obtain
$$
\|2BR(y,\zeta)\cdot y_\al-2BR(x,\zeta)\cdot x_\al+\frac{\rho_0}{\mu_0}(\grad P_2^{-1}(y)\cdot y_{\al}-\grad P_2^{-1}(x)\cdot x_{\al})\|_{L^2}\leq C\|z\|_{H^1},
$$
giving
$$
I_{2,3}\leq C\|z\|_{H^1}^2,
$$
as desired. Next we move to $I_3$. We split further to deal with $c_x-c_y$ considering $c_x-c_y=G_1+G_2$ where
$$
G_1=\frac{\al+\pi}{2\pi}\int \Big[\partial_\beta(Q^2_xBR_x)(\beta)\cdot\frac{x_\beta(\beta)}{|x_\beta(\beta)|^2}
-\partial_\beta(Q^2_yBR_y)(\beta)\cdot\frac{y_\beta(\beta)}{|y_\beta(\beta)|^2}\Big]d\beta,
$$
and
$$
G_2=-\int_{-\pi}^\al \Big[\partial_\beta(Q^2_xBR_x)(\beta)\cdot\frac{x_\beta(\beta)}{|x_\beta(\beta)|^2}
-\partial_\beta(Q^2_yBR_y)(\beta)\cdot\frac{y_\beta(\beta)}{|y_\beta(\beta)|^2}\Big]d\beta.
$$
Then we decompose further, to find $|G_1|\leq |G_{1,1}|+|G_{1,2}|+|G_{1,3}|+|G_{1,4}|+|G_{1,5}|$ where
$$
G_{1,1}=\int \partial_\al((Q^2_x-Q^2_y)BR_x)\cdot\frac{x_\al}{|x_\al|^2}d\al,
\quad
G_{1,2}=\int \partial_\al(Q^2_y) (BR_x-BR_y)\cdot\frac{x_\al}{|x_\al|^2}d\al,
$$
$$
G_{1,3}=\int Q^2_y\partial_\al(BR_x-BR_y)\cdot\frac{x_\al}{|x_\al|^2}d\al,
\quad
G_{1,4}=\int \partial_\al(Q^2_yBR_y)\cdot\frac{z_\al}{|x_\al|^2}d\al,
$$
$$
G_{1,5}=\int \partial_\beta(Q^2_yBR_y)\cdot y_\al
\Big(\frac{1}{|x_\al|^2}-\frac{1}{|y_\al|^2}\Big)d\al.
$$
Above we use $\al$ variables instead of $\beta$ for the sake of simplicity.
We can proceed as before to get
$$
|G_{1,1}|+|G_{1,2}|+|G_{1,4}|+|G_{1,5}|\leq C\| z\|_{H^1}.
$$
For the most delicate term we have to split further:
$
G_{1,3}= G_{1,3,1}+G_{1,3,2}+G_{1,3,3}+G_{1,3,4}+G_{1,3,5}+G_{1,3,6},
$
where
$$
G_{1,3,1}=\frac{1}{2\pi}\int Q^2_y\int\frac{x_-^\bot}{|x_-|^2}\cdot\frac{x_\al}{|x_\al|^2}\omega'_\al d\beta d\al,\,\,
G_{1,3,2}=\frac{1}{2\pi}\int Q^2_y\int\big[\frac{x_-^\bot}{|x_-|^2}-\frac{y_-^\bot}{|y_-|^2}\big]\cdot\frac{x_\al}{|x_\al|^2}\zeta'_\al d\beta d\al,
$$
$$
G_{1,3,3}=\frac{1}{2\pi}\int Q^2_y\int\frac{\da z_-^\bot}{|x_-|^2}\cdot\frac{x_\al}{|x_\al|^2}\gamma' d\beta d\al,\,\,
G_{1,3,4}=\frac{1}{2\pi}\int Q^2_y\int \da y_-^\bot \big[\frac{\gamma'}{|x_-|^2}-\frac{\zeta'}{|y_-|^2}\big]\cdot\frac{x_\al}{|x_\al|^2} d\beta d\al,
$$
$$
G_{1,3,5}=-\frac{1}{\pi}\int Q^2_y\int\frac{x_-^\bot}{|x_-|^4}\cdot\frac{x_\al}{|x_\al|^2}x_-\cdot \da z_-\zeta'
d\beta d\al,
$$
and
$$
G_{1,3,6}=-\frac{1}{\pi}\int Q^2_y\int\big[\frac{x_-^\bot}{|x_-|^4}x_-\cdot\da y_-\gamma'-\frac{y_-^\bot}{|y_-|^4}y_-\cdot\da y_-\zeta'\big]\cdot\frac{x_\al}{|x_\al|^2}
d\beta d\al.
$$
We estimate first the less singular terms, which can be controlled as before as follows:
$$
|G_{1,3,2}|+|G_{1,3,4}|+|G_{1,3,6}|\leq C \|z\|_{H^1}.
$$
One could rewrite $G_{1,3,1}$ as follows:
$$
G_{1,3,1}=\frac{1}{2\pi}\int Q^2_y\int\frac{x_-^\bot-x_\al^\bot\beta}{|x_-|^2}\cdot\frac{x_\al}{|x_\al|^2}\omega'_\al d\beta d\al,\quad
$$
to find a kernel of degree $0$ applied to $\omega_\al$. This yields
$$
|G_{1,3,1}|\leq C \|\omega\|_{L^2}\leq C \|z\|_{H^1}.
$$
Similarly
$$
G_{1,3,5}=-\frac{1}{\pi}\int Q^2_y\int\frac{x_-^\bot-x_\al^\bot\beta}{|x_-|^4}\cdot\frac{x_\al}{|x_\al|^2}x_-\cdot \da z_-\zeta'
d\beta d\al,
$$
and a kernel of order $-1$ applied to $\da z$ yields
$$
|G_{1,3,5}|\leq C \|z\|_{H^1}.
$$
It remains to deal with $G_{1,3,3}$ where we simply integrate by parts to obtain
$$
G_{1,3,3}=-\frac{1}{2\pi}\int \int z_-^\bot\cdot\da\Big(\frac{1}{|x_-|^2}Q^2_y\frac{x_\al}{|x_\al|^2}\gamma'\Big) d\beta d\al.
$$
We find as before
$$
|G_{1,3,3}|\leq C \|z\|_{H^1}.
$$
Since we are done with $G_1$ it remains to deal with $G_2$. Then, the same decomposition is going to work to control $G_2$ in the same manner than $G_1$, but for the analogous to the term $G_{1,3,3}$:
$$
G_{2,3,3}=\frac{1}{2\pi}\int_{-\pi}^{\al} Q^2_y(\beta)\int\frac{(\partial_{\beta} z(\beta)-\partial_{\beta} z(\beta-\xi)) ^\bot}{|x(\beta)-x(\beta-\xi)|^2}\cdot\frac{x_\beta(\beta)}{|x_\beta(\beta)|^2}\gamma(\beta-\xi) d\xi d\beta.
$$
We can not integrate by parts here as in $G_{1,3,3}$. We decompose further $G_{1,3,3}=G_{1,3,3}^1+G_{1,3,3}^2+G_{1,3,3}^3$
where
\begin{align*}
G^1_{2,3,3}=\frac{1}{2\pi}\int_{-\pi}^{\al} Q^2_y(\beta)\frac{x_\beta(\beta)}{|x_\beta(\beta)|^2}\cdot\int(\partial_{\beta}& z(\beta)-\partial_{\beta} z(\beta-\xi)) ^\bot\times\\
&\times
\Big[\frac{\gamma(\beta-\xi)}{|x(\beta)-x(\beta-\xi)|^2}-\frac{\gamma(\beta)}{|x_\beta(\beta)|^24\sin^2(\beta/2)}\Big] d\xi d\beta,
\end{align*}
\begin{align*}
G^2_{2,3,3}=&\frac{1}{2}\int_{-\pi}^{\al} Q^2_y(\beta)\gamma(\beta)\frac{x_\beta(\beta)}{|x_\beta(\beta)|^4}\cdot\Lambda(\partial_{\beta} z^{\bot})(\beta)d\beta\\
&-\frac{\al+\pi}{4\pi}\int Q^2_y(\beta)\gamma(\beta)\frac{x_\beta(\beta)}{|x_\beta(\beta)|^4}\cdot\Lambda(\partial_{\beta} z^{\bot})(\beta)d\beta,
\end{align*}
and
$$
G^3_{2,3,3}=\frac{\al+\pi}{2}\int Q^2_y(\beta)\gamma(\beta)\frac{x_\beta(\beta)}{|x_\beta(\beta)|^4}\cdot\Lambda(\partial_{\beta} z^{\bot})(\beta)d\beta.
$$
The fact that the kernel in $\xi$ has degree $-1$ allows us to get
$$
|G^1_{2,3,3}|\leq C\|z\|_{H^1}.
$$
Integrating by parts $\Lambda$ as a self-adjoint operator it is easy to obtain
$$
|G^3_{2,3,3}|\leq C\|z\|_{H^1}.
$$
All the bounds above for $c_x-c_y$ allow us to get
$$
I_3\leq C\|z\|^2_{H^1}+\int z\cdot x_\al G^2_{2,3,3}d\al.
$$
Above we integrate by parts to find
$$
\int z\cdot x_\al G^2_{2,3,3}d\al=I_{3,1}+I_{3,2}
$$
where
$$
I_{3,1}=\frac{1}{2}\int \Big(\int_{-\pi}^\al z(\beta)\cdot x_\beta(\beta)d\beta\Big)Q^2_y\gamma\frac{x_\al}{|x_\al|^4}\cdot\Lambda(\partial_{\al} z^{\bot})d\al.
$$
and
$$
I_{3,2}=-\frac{1}{4\pi}\int \int_{-\pi}^\al z(\beta)\cdot x_\beta(\beta)d\beta d\al \,\int Q^2_y(\beta)\gamma(\beta)\frac{x_\beta(\beta)}{|x_\beta(\beta)|^4}\cdot\Lambda(\partial_{\beta} z^{\bot})(\beta)d\beta
$$
As before, using that $\Lambda$ is self-adjoint it is easy to get
$$
I_{3,2}\leq C\|z\|^2_{H^1}.
$$
Similarly
$$
I_{3,1}=\frac{1}{2}\int \Lambda\Big(a\,Q^2_y\gamma\frac{x_\al}{|x_\al|^4}\Big)\cdot \partial_{\al} z^{\bot}d\al,
\quad
\mbox{where}
\quad
a(\al)=\int_{-\pi}^\al z(\beta)\cdot x_\beta(\beta)d\beta.
$$
The fact that $\Lambda=H(\da)$ allows us to find
$$
I_{3,2}\leq C\|z\|^2_{H^1},
$$
and finally
$$
I_{3}\leq C\|z\|^2_{H^1}.
$$
At this point it is easy to get
$$
I_{4}\leq C\|z\|^2_{H^1}.
$$
If one gathers above inequalities the following is obtained:
$$
\frac{d}{dt}\|z\|^2_{L^2}\leq C\|z\|^2_{H^1}.
$$
Next step is to analyzed
$$
\frac{d}{dt}\|z_{\al}\|^2_{L^2}=2\int \da z\cdot \da z_t d\g=I_5+I_6+I_7+I_8+I_9,
$$
where
$$
I_5=2\int z_\al\cdot \da\big(Q^2_x\big)(BR_x-BR_y)d\al,\quad
I_{6}=2\int z_\al\cdot Q^2_x\da(BR_x-BR_y)d\al
$$
$$
I_7=2\int z_\al\cdot \da\big((Q^2_x-Q^2_y)BR_y\big)d\al,\quad
I_8=2\int z_\al\cdot \da\big((c_x-c_y)x_\al\big) d\al,
$$
$$
I_9=2\int z_\al\cdot \da\big(c_y z_\al\big) d\al.
$$
It is easy to get
$$
I_{5}\leq C\|z\|^2_{H^1}.
$$
For $I_6$ we consider $I_{6}=I_{6,1}+I_{6,2}+I_{6,3}+I_{6,4}+I_{6,5}+I_{6,6}$ where
$$
I_{6,1}=\frac{1}{\pi}\int z_\al\cdot Q^2_x\int\frac{x_-^\bot}{|x_-|^2}\omega'_\al d\beta d\al,\quad
I_{6,2}=\frac{1}{\pi}\int  z_\al\cdot Q^2_x\int\big[\frac{x_-^\bot}{|x_-|^2}-\frac{y_-^\bot}{|y_-|^2}\big]\zeta'_\al d\beta d\al,
$$
$$
I_{6,3}=\frac{1}{\pi}\int z_\al\cdot Q^2_x\int\frac{\da z_-^\bot}{|x_-|^2}\gamma' d\beta d\al,\quad
I_{6,4}=\frac{1}{\pi}\int z_\al\cdot Q^2_x\int \da y_-^\bot \big[\frac{\gamma'}{|x_-|^2}-\frac{\zeta'}{|y_-|^2}\big]d\beta d\al,
$$
$$
I_{6,5}=-\frac{2}{\pi}\int  z_\al\cdot Q^2_x\int\frac{x_-^\bot}{|x_-|^4} x_-\cdot \da z_-\zeta'
d\beta d\al,
$$
and
$$
I_{6,6}=-\frac{2}{\pi}\int z_\al\cdot Q^2_x\int\big[\frac{x_-^\bot}{|x_-|^4}x_-\cdot\da y_-\gamma'-\frac{y_-^\bot}{|y_-|^4}y_-\cdot\da y_-\zeta'\big]
d\beta d\al.
$$
It is easy to get
$$
I_{6,2}+I_{6,4}+I_{6,6}\leq C\|z\|^2_{H^1}.
$$
We split further: $I_{6,3}=I_{6,3,1}+I_{6,3,2}$ to find
$$
I_{6,3,1}=\frac{1}{2\pi}\int z_\al\cdot \int\frac{\da z_-^\bot}{|x_-|^2}[Q^2_x\gamma'-(Q_x^{2})' \gamma] d\beta d\al,
$$
and
$$
I_{6,3,2}=\frac{1}{2\pi}\int z_\al\cdot \int\frac{\da z_-^\bot}{|x_-|^2}[Q^2_x\gamma'+(Q_x^{2})' \gamma] d\beta d\al.
$$
We find as before
$$
I_{6,3,1}\leq C\|z\|^2_{H^1}.
$$
Changing variables one could obtain
$$
I_{6,3,2}=\frac{1}{4\pi}\int\int \da z_-\cdot \frac{\da z_-^\bot}{|x_-|^2}[Q^2_x\gamma'+(Q_x^{2})' \gamma] d\beta d\al=0.
$$
We are done with $I_{6,3}$. The term $I_{6,5}$ is decomposed as follows
$$
I_{6,5,1}=-\frac{2}{\pi}\int  z_\al\cdot Q^2_x\int  \Big[\frac{\zeta' x_-^\bot}{|x_-|^4} x_- - \frac{\zeta x_\al^\bot}{|x_\al|^44\sin^2(\beta/2)} x_\al\Big]\cdot\da z_-
d\beta d\al,
$$
$$
I_{6,5,2}=-2\int  z_\al\cdot Q^2_x\zeta \frac{ x_\al^\bot}{|x_\al|^4} x_\al\cdot\Lambda(z_\al)
d\al.
$$
Therefore it yields
$$
I_{6,5,1}\leq C\|z\|^2_{H^1}.
$$
We use the fact that $\Lambda=H(\da)$ and the identity
$$
x_\al\cdot\da^2z=-x_\al\cdot\da^2 y=-z_\al\cdot\da^2y
$$
to rewrite
\begin{align*}
I_{6,5,2}=&2\int  z_\al\cdot Q^2_x\zeta \frac{ x_\al^\bot}{|x_\al|^4} [\Lambda(x_\al\cdot z_\al)-x_\al\cdot\Lambda(z_\al)]-2\int  z_\al\cdot Q^2_x\zeta \frac{x_\al^\bot}{|x_\al|^4} H(\da^2x\cdot z_\al)
d\al\\
&+2\int  z_\al\cdot Q^2_x\zeta \frac{x_\al^\bot}{|x_\al|^4} H(z_\al\cdot\da^2y )
d\al
\end{align*}
Since the commutator estimate for $\Lambda$ gives
$$
I_{6,5,2}\leq C\|z\|^2_{H^1},
$$
we are done with $I_{6,5}$. It remains to deal with $I_{6,1}$ where we have to find the Rayleigh-Taylor condition. We decompose further
$$
I_{6,1,1}=\frac{1}{\pi}\int z_\al\cdot Q^2_x\int\big[\frac{x_-^\bot}{|x_-|^2}-\frac{x_\al^\bot}{|x_\al|^22\tan(\beta/2)}\big]\omega'_\al d\beta d\al,\quad
I_{6,1,2}=\int Q^2_x z_\al\cdot\frac{x_\al^\bot}{|x_\al|^2} H(\omega_\al)d\al.
$$
As before
$$
I_{6,1,1}\leq C\|z\|^2_{H^1}.
$$
Next we will decompose $\omega_\al$, pointing out first the bounded terms, and dealing later with the unbounded. We take
$\omega_\al=G_3+G_4+G_5+G_6$
where
$$
G_3=-2BR_x\cdot \da^2x+2BR_y\cdot \da^2y,\quad
G_4=-2\da BR_x\cdot x_\al+2\da BR_y\cdot y_\al
$$
$$
G_5=-2\frac{\rho_0}{\mu_0}(\da(\grad P_2^{-1}(x))\cdot x_\al-\da(\grad P_2^{-1}(y))\cdot y_\al),\quad G_6=-2\frac{\rho_0}{\mu_0}(\grad P_2^{-1}(x)\cdot \da^2x-\grad P_2^{-1}(y)\cdot \da^2y).
$$
We split further $G_3=G_{3,1}+G_{3,2}$
$$
G_{3,1}=-2BR_x\cdot \da^2z,\quad G_{3,2}=2(BR_y-BR_x)\cdot \da^2y
$$
to obtain as before
$$
\|G_{3,2}\|_{L^2}\leq C\|z\|_{H^1}.
$$
The term $G_{3,1}$ is part of the unbounded characters. We continue by taking $G_4=G_{4,1}+G_{4,2}+G_{4,3}+G_{4,4}+G_{4,5}+G_{4,6}$
where
$$
G_{4,1}=-\frac1{\pi}\int\frac{\da z_-^\bot}{|x_-|^2}\gamma'd\beta \cdot x_\al,\quad
G_{4,2}=-\frac1{\pi}\int\da y_-^\bot\cdot\big[\frac{\gamma'}{|x_-|^2} x_\al-\frac{\zeta'}{|y_-|^2} y_\al\big]d\beta,
$$
$$
G_{4,3}=\frac2{\pi}\int\frac{x_-^\bot}{|x_-|^4}\cdot x_\al\, x_-\cdot \da z_-\gamma'd\beta ,\quad
G_{4,4}=\frac2{\pi}\int\big[\frac{x_-^\bot}{|x_-|^4}\cdot x_\al\, x_-\gamma'-\frac{y_-^\bot}{|y_-|^4}\cdot y_\al\, y_-\zeta'\big]\cdot \da y_-d\beta,
$$
$$
G_{4,5}=-\frac1{\pi}\int\frac{x_-^\bot}{|x_-|^2}\omega'_\al d\beta \cdot x_\al,\quad
G_{4,6}=-\frac1{\pi}\int \big[\frac{x_-^\bot}{|x_-|^2}\cdot x_\al-\frac{y_-^\bot}{|y_-|^2}\cdot y_\al\big]\zeta'_\al d\beta.
$$
Next, $G_{4,1}$ joints the unbounded terms and
$$
\|G_{4,2}\|_{L^2}+\|G_{4,4}\|_{L^2}+\|G_{4,6}\|_{L^2}\leq C\|z\|_{H^1}.
$$
It is possible to obtain a kernel of degree $-1$ applied to $\da z$ in $G_{4,3}$ as follows:
$$
G_{4,3}=\frac2{\pi}\int\frac{x_-^\bot-x^\bot_\al\beta}{|x_-|^4}\cdot x_\al\, x_-\cdot \da z_-\gamma'd\beta.
$$
Therefore
$$
\|G_{4,3}\|_{L^2}\leq C\|z\|_{H^1}.
$$
Since
$$
G_{4,5}=-\frac1{\pi}\int\frac{x_-^\bot-x_\al^\bot\beta}{|x_-|^2}\omega'_\al d\beta \cdot x_\al
$$
we obtain in an analogous way
$$
\|G_{4,5}\|_{L^2}\leq C\|\omega\|_{L^2}\leq C\|z\|_{H^1}.
$$
For $G_5$ it is easy to get
$$
\|G_{5}\|_{L^2}\leq C\|z\|_{H^1},
$$
but the term $G_{6}$ has to be decomposed as follows:
$$
G_{6,1}=-2\frac{\rho_0}{\mu_0}\grad P_2^{-1}(x)\cdot \da^2z,\quad G_{6,2}=-2\frac{\rho_0}{\mu_0}(\grad P_2^{-1}(x)-\grad P_2^{-1}(y))\cdot \da^2y.
$$
$G_{6,1}$ remains unbounded and
$$
\|G_{6,2}\|_{L^2}\leq C\|z\|_{H^1},
$$
easily. Thanks to all this decomposition we find
$$
I_{6,1,2}\leq C\|z\|^2_{H^1}+I_{6,1,2}^1+I_{6,1,2}^2+I_{6,1,2}^3,
$$
where
$$
I_{6,1,2}^1=\int Q^2_x z_\al\cdot\frac{x_\al^\bot}{|x_\al|^2} H(G_{3,1})d\al,\quad I_{6,1,2}^2=\int Q^2_x z_\al\cdot\frac{x_\al^\bot}{|x_\al|^2} H(G_{4,1})d\al,
$$
and
$$
I_{6,1,2}^3=\int Q^2_x z_\al\cdot\frac{x_\al^\bot}{|x_\al|^2} H(G_{6,1})d\al.
$$
In $I_{6,1,2}^1$ and $I_{6,1,2}^3$ the Rayleigh-Taylor condition will show up. For $I_{6,1,2}^2$ we consider the splitting
$$
I_{6,1,2}^{2,1}=\int H(Q^2_x z_\al\cdot\frac{x_\al^\bot}{|x_\al|^2})\frac1{\pi}\int\da z_-^\bot\big[\frac{\gamma'}{|x_-|^2}-\frac{\gamma}{|x_\al|^24\sin^2(\beta/2)}\big]d\beta \cdot x_\al  d\al,
$$
$$
I_{6,1,2}^{2,2}=\int H(Q^2_x z_\al\cdot\frac{x_\al^\bot}{|x_\al|^2})\frac{\gamma}{|x_\al|^2}\Lambda(\da z^\bot)\cdot x_\al  d\al,
$$
using that $H$ is skew-adjoint. First term satisfies
$$
I_{6,1,2}^{2,1}\leq C\|z\|^2_{H^1}.
$$
For the second one we use the commutator estimates to find
$$
I_{6,1,2}^{2,2}\leq C\|z\|^2_{H^1}+\int H(Q^2_x z_\al\cdot\frac{x_\al^\bot}{|x_\al|^2})\Lambda(\frac{\gamma}{|x_\al|^2}\da z^\bot\cdot x_\al)  d\al.
$$
The fact that $H^2=-I$ yields
$$
I_{6,1,2}^{2,2}\leq C\|z\|^2_{H^1}+\int Q^2_x z_\al\cdot\frac{x_\al^\bot}{|x_\al|^2}\da(\frac{\gamma}{|x_\al|^2}\da z^\bot\cdot x_\al)  d\al.
$$
In the integral above we expand the derivative, to find out that it is possible to integrate by parts in $\da(z_\al\cdot x_\al^\bot)$. This yields
$$
I_{6,1,2}^{2,2}\leq C\|z\|^2_{H^1},\quad \mbox{and therefore}\quad I_{6,1,2}^{2}\leq C\|z\|^2_{H^1}.
$$
Next we consider
$$
I_{6,1,2}^1=-\int Q^2_x z_\al\cdot\frac{x_\al^\bot}{|x_\al|^2} H(2BR_x\cdot \da^2z)d\al,
$$
for which we use the commutator for the Hilbert transform
$$
\|H(g\da f)-gH(\da f)\|_{L^2}\leq C\|g\|_{C^{1,\frac13}}\|f\|_{L^2},
$$
to find
$$
I_{6,1,2}^1\leq -\frac{2}{|x_\al|^2}\int Q^2_x z_\al\cdot x_\al^\bot BR_x\cdot H(\da^2z)d\al,
$$
Next we split above integral by components:
\begin{equation}\label{BRBR1}
I_{6,1,2}^{1,1}=\frac{2}{|x_\al|^2}\int Q^2_x \da z_1 \da x_2 BR_{x1}\cdot H(\da^2z_1)d\al,
\end{equation}
$$
I_{6,1,2}^{1,2}=\frac{2}{|x_\al|^2}\int Q^2_x \da z_1 \da x_2 BR_{x2}\cdot H(\da^2z_2)d\al,
$$
$$
I_{6,1,2}^{1,3}=-\frac{2}{|x_\al|^2}\int Q^2_x \da z_2 \da x_1 BR_{x1}\cdot H(\da^2z_1)d\al,
$$
\begin{equation}\label{BRBR4}
I_{6,1,2}^{1,4}=-\frac{2}{|x_\al|^2}\int Q^2_x \da z_2 \da x_1 BR_{x2}\cdot H(\da^2z_2)d\al.
\end{equation}
The commutator for the Hilbert transform allows us to obtain
$$
I_{6,1,2}^{1,2}\leq C\|z\|^2_{H^1}+\frac{2}{|x_\al|^2}\int Q^2_x \da z_1  BR_{x2}\cdot H(\da x_2\da^2z_2)d\al,
$$
and together with identity
$$
\da x_2\da^2z_2=-\da x_1\da^2z_1-\da z\cdot \da^2 y
$$
provides
$$
I_{6,1,2}^{1,2}\leq C\|z\|^2_{H^1}-\frac{2}{|x_\al|^2}\int Q^2_x \da z_1  BR_{x2}\cdot H(\da x_1\da^2z_1)d\al.
$$
The commutator estimate yields
\begin{equation}\label{BRBR2}
I_{6,1,2}^{1,2}\leq C\|z\|^2_{H^1}-\frac{2}{|x_\al|^2}\int Q^2_x \da z_1  \da x_1 BR_{x2}\cdot H(\da^2z_1)d\al.
\end{equation}
In a similar manner we find
\begin{equation}\label{BRBR3}
I_{6,1,2}^{1,3}\leq C\|z\|^2_{H^1}+\frac{2}{|x_\al|^2}\int Q^2_x \da z_2  \da x_2 BR_{x1}\cdot H(\da^2z_2)d\al.
\end{equation}
Adding \eqref{BRBR1}, \eqref{BRBR2}, \eqref{BRBR3} and \eqref{BRBR4} we find
\begin{equation}\label{BRF}
I_{6,1,2}^1\leq C\|z\|^2_{H^1}-\frac{2}{|x_\al|^2}\int Q^2_x BR_x \cdot  x_\al^\bot  z_\al \Lambda(z_\al)d\al.
\end{equation}
Next
$$
I_{6,1,2}^3=-2\frac{\rho_0}{\mu_0}\int Q^2_x z_\al\cdot\frac{x_\al^\bot}{|x_\al|^2} H(\grad P_2^{-1}(x)\cdot \da^2z)d\al,
$$
and a decomposition in components as before provides
\begin{equation}\label{BRP}
I_{6,1,2}^3\leq C\|z\|^2_{H^1}-2\frac{\rho_0}{\mu_0|x_\al|^2}\int Q^2_x \grad P_2^{-1}(x)\cdot x_\al^\bot z_\al\cdot \Lambda(z_\al)d\al.
\end{equation}
Adding \eqref{BRF} and \eqref{BRP} we find
\begin{equation*}
I_{6,1,2}\leq C\|z\|^2_{H^1}+I_{6,1,2}^{1}+I_{6,1,2}^{3}\leq C\|z\|^2_{H^1}-2\frac{1}{\mu_0|x_\al|^2}\int Q^2_x \sigma_x z_\al\cdot \Lambda(z_\al)d\al.
\end{equation*}
The positivity of the Rayleigh-Taylor condition for the curve $x$ gives
$$
I_{6,1,2}\leq C\|z\|^2_{H^1},\quad I_{6,1}\leq C\|z\|^2_{H^1},\quad\mbox{and finally}\quad I_{6}\leq C\|z\|^2_{H^1}.
$$
We find easily
$$
I_{7}\leq C\|z\|^2_{H^1}.
$$
For $I_8$ we consider
$$
I_{8}=2\int z_\al\cdot (c_x-c_y)\da^2 x d\al+2\int z_\al\cdot x_\al \da(c_x-c_y)  d\al,
$$
and integrate by parts to find
$$
I_{8}=-2\int \da^2 z\cdot x_\al (c_x-c_y)  d\al.
$$
The fact that
$$
I_8=2\int \da^2 y\cdot z_\al (c_x-c_y)  d\al,
$$
allows us to deal with $I_{8}$ as for $I_3$ to get
$$
I_{8}\leq C\|z\|_{H^1}^2.
$$
Finally, integration by parts provides
$$
I_{9}=\int |z_\al|^2\da c_y  d\al\leq C\|z\|_{H^1}^2.
$$
\section{Applying perturbative argument and concluding the proof}
Finally, we will applied a perturbative argument to conclude the proof of theorem \ref{splashvacuum}. Consider
a curve $z^l(\al)$ as in section 2 and $P(z^l(\al))$ an initial datum for $P$(Muskat). Then we get a solution $\widetilde{z}^l(\al,t)\in C([0,T],H^3)$ given by using section 3. Next we consider a perturbation of $z^l(\al)$, the curve $z_0(\al)$, for which the Rayleigh-Taylor and chord-arc conditions holds. Furthermore, the velocity given by Darcy's law for $z_0(\al)$ shows that two different branches of the interface are going to approach as time goes forward. Next we take $P(z_0(\al))$ as an initial datum for $P$(Muskat), getting a solution $\widetilde{z}(\al,t)\in C([0,T],H^3)$. The stability result in section 5 gives
$$
\|\widetilde{z}-\widetilde{z}^l\|_{H^1}(t)\leq C(\sup_{[0,T]}E_3(\widetilde{z},t)+\sup_{[0,T]}E_3(\widetilde{z}^l,t))\|P(z_0(\al))-P(z^l(\al))\|_{H^1}
$$
and therefore
$$
\|\widetilde{z}-\widetilde{z}^l\|_{H^1}(t)\leq C(\sup_{[0,T]}E_3(\widetilde{z},t)+\sup_{[0,T]}E_3(\widetilde{z}^l,t))\|z_0(\al)-z^l(\al)\|_{H^1}.
$$
Here we point out that the time of existence in section 3 is independent of the smallness of $\|z_0(\al)-z^l(\al)\|_{H^1}$. Since the transformation $P^{-1}$ is well define for $\tilde{z}$ and the fact that $z^l$ self-intersect at a point allows us to conclude that in the evolution of $z=P^{-1}(\widetilde{z})$ there exists a finite time such that $z$ has to break down with a splash singularity.

\section{A remark on the family of splash singularities }\label{generales}
The  scenario in section \ref{initialdata} is the simplest one to obtain a splash singularity. However the one-phase Muskat problem can develop this kind of point-wise collapse for more geometries. In order to check that we proceed as follows. Let $z(\alpha)$ be a splash curve with $\alpha_1\neq \alpha_2$ such that $z(\alpha_1)=z(\alpha_2)$ and $|\pa_\alpha z (\alpha)|>0$ for every $\alpha$. To consider a different situation than in section 2, we also assume that $\pa_\alpha z_1(\alpha_1)\neq 0$. We make the following distinction between $\alpha_1$ and $\alpha_2$:
There exist a neighborhood $U_{\alpha_1}$ of $\alpha_1$ and a neighborhood $U_{\alpha_2}$ of $\alpha_2$ such that, if $z_1(\beta_1)=z_1(\beta_2)$ for $\beta_1\in U_{\alpha_1}$ and $\beta_2\in U_{\alpha_2}$, then $z_2(\beta_1)\leq z_2(\beta_2)$. Roughly speaking, we just mean that $z(\alpha_2)$ is the upper splash point and $z(\alpha_1)$ is the lower splash point.

Let us analyze the normal velocity at $\alpha_2$ and $\alpha_1$. For $\alpha_2$ we have
$$\mu_0u(\alpha_2)\cdot n(\alpha_2)=-\grad p(z(\alpha_2))\cdot n(\al_2)-\rho_0 \frac{\pa_\alpha z_1(\alpha_2)}{|\pa_\alpha z(\alpha_2)|},$$
where $n(\al)=\pa_\alpha^{\bot} z(\alpha)/|\pa_\alpha z(\alpha)|$.
As in section \ref{initialdata}, Hopf's lemma yields $$-\grad p(z(\alpha_2))\cdot n(\alpha_2)>0$$ and we consider $$\frac{\pa_\alpha z_1(\alpha_2)}{|\pa_\alpha z(\alpha_2)|}< 0.$$
On the other hand, we have
$$\mu_0u(\alpha_1)\cdot n(\alpha_1)=-\grad p(z(\alpha_1))\cdot n(\alpha_1)-\rho_0 \frac{\pa_\alpha z_1(\alpha_1)}{|\pa_\alpha z(\alpha_1)|},$$
with
$$-\grad p(z(\alpha_1))\cdot n(\alpha_1)>0,$$
and
$$\frac{\pa_\alpha z_1(\alpha_1)}{|\pa_\alpha z(\alpha_1)|}> 0.$$
Then the sign of $\pa_\alpha z(\alpha_1)$ is bad for our purpose. However we can notice that
$$\frac{\pa_\alpha z_1(\alpha_1)}{|\pa_\alpha z(\alpha_1)|}=- \frac{\pa_\alpha z_1(\alpha_2)}{|\pa_\alpha z(\alpha_2)|},$$
and therefore
$$u(\alpha_2)\cdot n(\alpha_2)> -u(\alpha_1)\cdot n(\alpha_1).$$
The last inequality is enough to show that the velocity separates the splash points backward in time. Unfortunately this is not enough to assure that we can produce a splash singularity by using the previous analysis. It is possible to find $u(\alpha_1)\cdot n(\alpha_1)$ negative. Then the solution would cross the branch of $P$ backward in time. This is a mere technical problem that we can solve as follow.

Let's define a velocity $$\underline{v}(x_1,x_2,t)=v(x_1,x_2-\frac{\rho_0}{\mu_0} t,t)+(0,\frac{\rho_0}{\mu_0}),$$
a density $$\underline{\rho}(x_1,x_2,t)=\rho(x_1,x_2-\frac{\rho_0}{\mu_0} t,t),$$
and a viscosity
$$\underline{\mu}(x_1,x_2,t)=\mu(x_1,x_2-\frac{\rho_0}{\mu_0}t,t).$$ Therefore
\begin{align*}
\pa_t \underline{\rho}(x_1,x_2,t)=&(\pa_t\rho)(x_1,x_2-\frac{\rho_0}{\mu_0}t,t)-\frac{\rho_0}{\mu_0}(\pa_{x_2}\rho)(x_1,x_2-\frac{\rho_0}{\mu_0} t,t)\\
=& - v(x_1,x_2-\frac{\rho_0}{\mu_0},t)\cdot (\nabla\rho) (x_1,x_2-\frac{\rho_0}{\mu_0} t,t)-\frac{\rho_0}{\mu_0}(\pa_{x_2}\rho)(x_1,x_2-\frac{\rho_0}{\mu_0} t,t)\\
=& - v(x_1,x_2-\frac{\rho_0}{\mu_0},t)\cdot \nabla\underline{\rho} (x_1,x_2,t)-\frac{\rho_0}{\mu_0}\pa_{x_2}\underline{\rho}(x_1,x_2,t)\\
=&  - \left(v(x_1,x_2-\frac{\rho_0}{\mu_0}t,t)+(0,\frac{\rho_0}{\mu_0})\right)\cdot \nabla\underline{\rho}(x_1,x_2,t).
\end{align*}
Thus we have that $\underline{\rho}$ satisfies
\begin{align*}
\pa_t\underline{\rho}+ \underline{v}\cdot \nabla \underline{\rho}=0,
\end{align*}
and in a similar manner it is easy to get
$$
\pa_t\underline{\mu}+ \underline{v}\cdot \nabla \underline{\mu}=0.
$$
On the other hand
\begin{align*}
\underline{\mu}\,\underline{v}=-\nabla\underline{p},
\end{align*}
where we consider $$\underline{p}(x_1,x_2,t)=p(x_1,x_2-\frac{\rho_0}{\mu_0}t,t).$$
Then,  by using $\underline{v}$, $\underline{\rho}$, $\underline{\mu}$ and $\underline{p}$,  we can write our Muskat problem as the system
\begin{align*}
\pa_t\underline{\rho}+ \underline{v}\cdot \nabla \underline{\rho}=&0,\\
\pa_t\underline{\mu}+ \underline{v}\cdot \nabla \underline{\mu}=&0,\\
\underline{\mu}\,\underline{v}=&-\nabla \underline{p},\\
\nabla \cdot \underline{v}=&0,
\end{align*}
with the boundary condition $$\lim_{x_2\to -\infty} \left(\underline{v}(x_1,x_2,t)-(0,\frac{\rho_0}{\mu_0})\right)=0.$$

In this new system we find the following: If $z(\alpha)$ is a splash curve such that $z(\alpha_1)=z(\alpha_2)$ then
\begin{align*}
\mu_0\underline{v}(z(\alpha_1))\cdot n(\alpha_1)=& -\grad\underline{p}(z(\alpha_1))\cdot n(\alpha_1),\\
\mu_0\underline{v}(z(\alpha_2))\cdot n(\alpha_2)=& -\grad\underline{p}(z(\alpha_2))\cdot n(\alpha_2),
\end{align*}
and again we can invoke Hopf's lemma to obtain that
\begin{align*}
-\grad\underline{p}(z(\alpha_1))\cdot n(\alpha_1)>0,\quad -\grad\underline{p}(z(\alpha_2))\cdot n(\alpha_2)>0.
\end{align*}
Then, the velocity separates the splash point and  $\underline{u}(\alpha_1)\cdot n(\alpha_1)n(\alpha_1)$ points in the opposite direction to $\underline{u}(\alpha_2)\cdot n(\alpha_2)n(\alpha_2)$. Therefore we can carry out the same analysis we did for the simpler case of section \ref{initialdata}.


\subsection*{{\bf Acknowledgements}}

\smallskip

AC, DC and FG were partially supported by the grant MTM2011-26696 (Spain). AC was partially supported by the ERC grant 307179-GFTIPFD. CF was partially supported by NSF and ONR grants DMS 09-01040 and N00014-08-1-0678. FG acknowledges support from the Ram\'on y Cajal program.


\begin{tabular}{ll}
\textbf{Angel Castro} & \textbf{Diego C\'ordoba} \\
{\small Instituto de Ciencias Matem\'aticas} & {\small Instituto de Ciencias Matem\'aticas}\\
{\small Universidad Aut\'onoma de Madrid} & {\small Consejo Superior de Investigaciones Cient\'ificas}\\
{\small C/ Nicolas Cabrera, 13-15, 28049 Madrid, Spain} & {\small C/ Nicolas Cabrera, 13-15, 28049 Madrid, Spain}\\
{\small Email: angel\underline{  }castro@icmat.es} & {\small Email: dcg@icmat.es}\\
   & \\
\textbf{Charles Fefferman} & \textbf{Francisco Gancedo}\\
{\small Department of Mathematics} & {\small Departamento de An\'{a}lisis Matem\'{a}tico $\&$ IMUS}\\
{\small Princeton University} & {\small Universidad de Sevilla}\\
{\small 1102 Fine Hall, Washington Rd, } & {\small C/ Tarfia s/n, Campus Reina Mercedes,}\\
{\small Princeton, NJ 08544, USA} & {\small 41012 Sevilla, Spain}\\
 {\small Email: cf@math.princeton.edu} & {\small Email: fgancedo@us.es}
  \\

\end{tabular}

\end{document}